\documentclass[conference]{ieeeconf}
\IEEEoverridecommandlockouts
\usepackage{cite}
\usepackage{amsmath,amssymb,amsfonts}
\usepackage{graphicx}
\usepackage{textcomp}
\usepackage{xcolor}
\usepackage{algorithm}
\usepackage{algpseudocode}
\usepackage{tabularray}
\usepackage{diagbox} 
\usepackage{array}   
\usepackage[colorlinks,
linkcolor=black,
anchorcolor=black,
citecolor=black]{hyperref}
\UseTblrLibrary{diagbox}

\newcommand{\ie}{\emph{i.e.}}

\renewcommand{\leq}{\leqslant}
\renewcommand{\geq}{\geqslant}
\renewcommand{\epsilon}{\varepsilon}
\renewcommand{\phi}{\varphi}
\renewcommand{\epsilon}{\varepsilon}

\begin{document}

\title{Efficient Simulation of Singularly Perturbed Systems Using a Stabilized Multirate Explicit Scheme 
\thanks{This work was partially supported by the Swedish Research Council under contract number 2023-05170. The authors are with the Division of Decision and Control Systems, KTH Royal Institute of Technology. Emails: yibos@kth.se (Y. S.), crro@kth.se (C. R.).}}

\author{\authorblockN{Yibo Shi and Cristian R. Rojas}}

\maketitle

\begin{abstract}
Singularly perturbed systems (SPSs) are prevalent in engineering applications, where numerically solving their initial value problems (IVPs) is challenging due to stiffness arising from multiple time scales. Classical explicit methods require impractically small time steps for stability, while implicit methods developed for SPSs are computationally intensive and less efficient for strongly nonlinear systems. This paper introduces a Stabilized Multirate Explicit Scheme (SMES) that stabilizes classical explicit methods without the need for small time steps or implicit formulations. By employing a multirate approach with variable time steps, SMES allows the fast dynamics to rapidly converge to their equilibrium manifold while slow dynamics evolve with larger steps. Analysis shows that SMES achieves numerical stability with significantly reduced computational effort and controlled error. Its effectiveness is illustrated with a numerical example.
\end{abstract}

\section{Introduction}
Singularly perturbed systems (SPSs) are prevalent in various engineering fields, including electrical and robotic systems, chemical reactions, and biological processes~\cite{kokotovic1999singular, khalil2002nonlinear}. These systems are characterized by multiple time scales due to small parameters multiplying the derivatives of certain state variables, leading to a separation between fast and slow dynamics. The coexistence of different time scales complicates the application of standard system analysis and control methods because the fast dynamics can introduce numerical stiffness and rapid transients~\cite{tikhonov1952systems, omalley1991singular}. Consequently, specialized approaches are necessary to handle the complexities introduced by these small perturbation parameters.

Extensive research has been devoted to various aspects of SPSs, such as model reduction techniques, composite control strategies, and asymptotic methods to simplify the analysis and controller design for such systems~\cite{kokotovic1999singular, kokotovic1984applications}. Within these problems, numerically solving initial value problems (IVPs) for SPSs plays a vital role. Numerical solutions of IVPs enable the exploration of system behavior, validation of theoretical results, and development of control algorithms when analytical solutions are intractable \cite{verhulst2005methods}. 

Numerical methods for solving IVPs of SPSs have been widely studied. Commonly used stiff solvers, such as MATLAB's \texttt{ODE15s} and \texttt{ODE23s}, are designed to handle stiff differential equations \cite{shampine1997matlab}. \texttt{ODE15s} uses Gear’s method, and \texttt{ODE23s} is based on a modified second-order Rosenbrock formula; both are fully implicit methods. These solvers provide numerical stability when simulating SPSs by allowing larger time steps without sacrificing significant accuracy. To address the computational burden of fully implicit methods, implicit-explicit (IMEX) schemes for SPSs have been developed~\cite{ascher1997implicit,akrivis2013implicit}. IMEX methods split the system into fast and slow components, treating the fast dynamics implicitly and the slow dynamics explicitly. This approach leverages the stability of implicit methods for the fast dynamics while retaining the computational efficiency of explicit methods for the slow dynamics \cite{boscarino2007error, pareschi2005implicit}.

However, the aforementioned methods for SPSs rely entirely or partially on implicit numerical methods. This is because the fast dynamics introduced by small time constants require indispensable numerical stability when using step sizes much larger than the time constants~\cite{cellier2006continuous}. Common explicit numerical methods, such as Forward Euler and the explicit Runge-Kutta methods, often require prohibitively small time steps to stabilize the fast dynamics, leading to increased computational cost and inefficiency~\cite{hairer1991solving}. Meanwhile, implicit methods, although more stable for SPSs, involve solving nonlinear equations at each iteration, which can be computationally intensive and may become inefficient for large-scale or highly nonlinear systems~\cite{ascher1998computer}. IMEX methods can take advantage of the efficiency of explicit methods and the stability of implicit methods; however, they require an explicit separation of the fast and slow dynamics, which is not always simple or trivial. 

The limitations of current numerical methods regarding efficiency, stability, and simplicity highlight the need for specialized techniques tailored to SPSs. Specifically, there is a demand for numerical schemes that can efficiently simulate these systems while providing stability guarantees and error control, without the computational burden of implicit methods or the requirement of time-scale separation as in IMEX methods. In this paper, we address this challenge by proposing an efficient Stabilized Multirate Explicit Scheme (SMES) designed for SPSs.

An important observation is that SPSs can usually be related to differential-algebraic equations (DAEs) when the perturbation parameter approaches zero~\cite{hairer1991solving}. In this limit, the fast dynamics effectively become algebraic constraints, and the system exhibits characteristics similar to DAEs. Numerical methods for solving DAEs often employ variable-step or multirate strategies to efficiently handle the coupled differential and algebraic components without resorting to excessively small time steps~\cite{brenan1996numerical}. Inspired by these techniques, SMES incorporates a multirate approach with variable time steps, which stabilizes classical explicit numerical methods to handle the fast dynamics. This approach maintains the simplicity and ease of implementation of explicit methods while enhancing computational efficiency and providing numerical stability guarantees. Moreover, since SMES is based on explicit schemes, it is better suited for systems with strong nonlinearity compared to implicit methods.

The remainder of this paper is organized as follows: In Section~\ref{sec:preliminary}, we provide the relevant background on SPSs. In Section~\ref{sec:FEM}, we present the limitation of classical explicit methods on SPSs using the Forward Euler method. In Section~\ref{sec:SMFE}, we introduce the SMES based on the Forward Euler method, outlining its theoretical foundation and implementation details. In Sections~\ref{sec:stability_analysis} and~\ref{sec:error_analysis}, we perform stability and error analyses of the proposed method. In Section~\ref{sec:numerical}, we present numerical illustrations to demonstrate the efficiency of our method on SPSs. Finally, Section~\ref{sec:conclusion} concludes the paper.

\textbf{Notation.} \textit{If $f$ and $g$ are real-valued functions defined on a subset of $\mathbb{R}_0^+$ including a neighborhood of $0$, we write $f = O(g)$ to mean that there are positive constants $K$ and $\delta$ such that
\begin{displaymath}
\vert f(h)\vert \le K \vert g(h)\vert, \quad \text{whenever}\quad 0 < h < \delta .
\end{displaymath}
If $f(h) = O(h^p)$ we say that ``$f$ is of order $h^p$''.}


\section{Preliminary}\label{sec:preliminary}

SPSs are dynamical systems characterized by multiple time scales. These systems are generally described by the state-space form \cite{levinson1950perturbations}
\begin{equation}\label{eq:generalss}
    \dot{X} = F(X, u, \epsilon, t), \quad X(t_0) = X^0,
\end{equation}
where $X \in \mathbb{R}^{n+m}$ is the state vector, $u$ is the input vector, $F$ is a sufficiently smooth function, and $\epsilon$ is a small singular perturbation parameter ($0 < \varepsilon \ll 1$). An analysis of \eqref{eq:generalss} is often based on the time-scale separated form~\cite{kokotovic1999singular}
\begin{equation}\label{eq:ctss}
    \begin{aligned}
        \dot{x} &= f(x, z, u, \epsilon, t), \quad x(t_0) = x^0, \\
        \varepsilon \dot{z} &= g(x, z, u, \epsilon, t), \quad z(t_0) = z^0,
    \end{aligned}
\end{equation}
where $x \in \mathbb{R}^n$ and $z \in \mathbb{R}^m$ are the separated slow and fast state variables. In this paper, we assume that the system in \eqref{eq:generalss} or \eqref{eq:ctss} is asymptotically stable.

For convenience, we will mainly consider SPSs in \textit{standard form}. When $\epsilon = 0$, the dimension of the state space in \eqref{eq:ctss} reduces from $n+m$ to $n$ because the differential equation of $\dot{z}$ degenerates into the algebraic equation
\begin{equation}\label{eq:algebric}
    0 = g(\Bar{x}, \Bar{z}, \Bar{u}, 0, t),
\end{equation}
where the bar is used to indicate that the variables belong to a system with $\epsilon = 0$.
Eq.~\eqref{eq:ctss} is said to be in standard form if and only if the following assumption concerning \eqref{eq:algebric} is satisfied:

\textbf{Assumption 1.} \cite{kokotovic1984applications} \textit{In a domain of interest, \eqref{eq:algebric} has $k\geq1$ distinct real roots
\begin{equation}\label{eq:assumption1}
    \Bar{z} = \Bar{\phi}_i(\Bar{x}, \Bar{u}, t), \quad i= 1, 2, \dots, k.
\end{equation} }

This assumption ensures a well-defined $n$-dimensional \textit{reduced model} of the full system in \eqref{eq:ctss} corresponding to each root, described as the differential-algebraic equation (DAE)
\begin{equation}\label{eq:reduced}
\begin{aligned}
    \dot{\Bar{x}} &= f(\Bar{x}, \Bar{z}, \Bar{u}, t), \quad x(t_0) = x^0, \\
    \Bar{z} &= \Bar{\phi}_{i}(\Bar{x}, \Bar{u}, t).
\end{aligned}
\end{equation}

Eq.~\eqref{eq:reduced} is also called a \textit{quasi-steady-state form} of \eqref{eq:ctss} because the fast dynamics in $z$, whose rate $\dot{z} = g/\epsilon$ can be large when $\epsilon$ is small, may rapidly converge to the equilibrium manifold $\mathcal{M}_0\colon \Bar{z} = \Bar{\phi}_{i}(\Bar{x}, \Bar{u}, t)$.

Typically, the reduced model represents the slow (average) phenomena that are dominant in most applications \cite{kokotovic1999singular}. Moreover, the reduced model does not impose strict stability requirements on the numerical methods used. However, since the quasi-steady-state $[\Bar{x}, \Bar{z}]$ is not free to start from $[x^0, z^0]$ as in \eqref{eq:ctss}, a discrepancy arises between the responses of the reduced and full models when $[x^0, z^0] \notin \mathcal{M}_0$, which lead to a fast transient. Therefore, to accurately capture the dynamics of the full state, it is necessary to study numerical methods applied to the full model in \eqref{eq:ctss}.

\textbf{Remark 1.} \textit{In this paper, we consider the more general case where the 
system may not be in standard form, or
may not have an explicit separation of fast and slow dynamics as in \eqref{eq:ctss}, which is not always straightforward to determine. Nevertheless, to simplify the subsequent analysis, we consider the standard form of singularly perturbed systems that satisfy \textbf{Assumption 1}.}

\section{Forward Euler Method}\label{sec:FEM}

In this section, we analyze the numerical stability of the Forward Euler Method (FEM) \cite{hairer1993solving}, as an example of a classical explicit method for SPSs, to demonstrate their inefficiency. Applying the classical FEM to solve the IVP in \eqref{eq:generalss} with step size $\Delta$ yields the discrete-time (DT) model 
\begin{equation}\label{eq:forwardeuler}
    X_{k+1} = X_k + \Delta F(X_k, u_k, \epsilon), \quad X_0 = X^0.
\end{equation}

To illustrate the effect of the small parameter $\epsilon$ on the numerical stability of FEM, we consider the following linear time-invariant autonomous model as a `test equation' \cite{kokotovic1999singular}:
\begin{equation}\label{eq:test1}
    \begin{aligned}
        \dot{x}(t) &= A_{11} x(t) + A_{12} z(t), \quad x(t_0) = x^0, \\
        \epsilon \dot{z}(t) &= A_{21} x(t) + A_{22} z(t), \quad z(t_0) = z^0,
    \end{aligned}
\end{equation}
where $x(t) \in \mathbb{R}^n$, $z(t) \in \mathbb{R}^m$, and the time variable $t$ is explicitly written out to facilitate subsequent analysis. If $A_{22}$ is nonsingular, there exists a matrix $L(\epsilon) \in \mathbb{R}^{m\times n}$ satisfying
\begin{align*}
    L(\varepsilon)&=A_{22}^{-1} A_{21}+\varepsilon A_{22}^{-2} A_{21} A_0+O\left(\varepsilon^2\right),
\end{align*}
where $A_0=A_{11}-A_{12} A_{22}^{-1} A_{21}$. This matrix $L(\epsilon)$ leads to the change of variables
\begin{equation}
    \eta(t) = z(t) + L(\epsilon)x(t).
\end{equation}

The system in \eqref{eq:test1} can then be rewritten in block-triangular form \cite{kokotovic1984applications}, to separate the slow and fast dynamics, as
\begin{equation}\label{eq:testtri}
    \left[\begin{array}{l}
\dot{x}(t) \\
\dot{\eta}(t)
\end{array}\right]=A \left[\begin{array}{l}
x(t) \\
\eta(t)
\end{array}\right],
\end{equation}
with
\begin{equation*}
A = \left[\begin{array}{cc}
A_{11}-A_{12} L(\epsilon) & A_{12} \\
0 & \epsilon^{-1}A_{22}+ L(\epsilon) A_{12}
\end{array}\right].
\end{equation*}

The system in \eqref{eq:test1} has been partially decoupled into a slow subsystem involving $x(t)$ and a fast subsystem involving $\eta(t)$. Correspondingly, the eigenvalues $\lambda$ of matrix $A$ can be separated into $\lambda_\text{slow}$ for the slow subsystem and $\lambda_\text{fast}$ respectively for the fast subsystem, where
\begin{align*}
    & \det(\lambda_\text{slow} I - A_{11}-A_{12} L(\epsilon)) = 0,\\
    & \det(\lambda_\text{fast} I - \epsilon^{-1}A_{22}+ L(\epsilon) A_{12}) = 0,
\end{align*}
where it can be observed that $\lambda_\text{slow} = O(1)$ and $\lambda_\text{fast} = O(1 / \epsilon)$. Applying the FEM to the system in block-triangular form leads to the following DT model:
\begin{align*}
\left[\begin{array}{l}
x(t + \Delta) \\
\eta(t + \Delta)
\end{array}\right] \approx \left( I + \Delta A \right) \left[\begin{array}{l}
x(t) \\
\eta(t)
\end{array}\right].
\end{align*}

For the numerical stability of the FEM, it is necessary and sufficient that $|1 + \Delta \lambda| < 1$ for every eigenvalue $\lambda$ of $A$. Moreover, since the system is asymptotically stable (i.e., $\mathrm{Re}(\lambda)\leq 0$ for each $\lambda$) and $\mathrm{Re}(\lambda_\text{fast}) \ll \mathrm{Re}(\lambda_\text{slow})$ for each $\lambda_\text{fast}$ and $\lambda_\text{slow}$, the maximum value of $\Delta$ ensuring stability is governed by the eigenvalue $\hat{\lambda}_\text{fast}$, where
\begin{align}\label{eq:maxlamdafast}
    \hat{\lambda}_\text{fast} = \mathop{\arg\min}\limits_{\lambda_{\text{fast}, i}}\{\mathrm{Re}(\lambda_{\text{fast}, i})\colon i = 1,2,\dots,m\}.
\end{align}

Therefore, the necessary and sufficient condition for numerical stability becomes $|1 + \hat{\lambda}_\text{fast} \Delta| < 1$, which leads to
\begin{equation}
    \Delta < -\frac{2}{\hat{\lambda}_\text{fast}} = O(\epsilon).
\end{equation}

Thus, the step size $\Delta$ is constrained by the smallest time scale associated with the fast dynamics. This implies that as $\varepsilon \rightarrow 0$, the step size $\Delta$ must decrease significantly. Consequently, a large number of time steps are required to cover the desired simulation horizon $T$, i.e., $T / \varepsilon = O(1 / \varepsilon)$ function evaluations, leading to increased computational cost and simulation time.

\section{Stabilized Multirate Forward Euler Method}\label{sec:SMFE}

Given the limitations of classical explicit numerical methods, we introduce a SMES that handles singular perturbations without requiring prohibitively small step sizes, thus allowing for faster and stable simulations of SPSs using explicit methods. In this section, we illustrate how SMES works by augmenting the FEM with the proposed scheme, resulting in a stabilized multirate Forward Euler (SMFE) method.

Instead of using a fixed step size $\Delta$ as in the classical FEM, we apply a multirate discretization within one step of length $\Delta$. First, we perform $N$ iterations of \eqref{eq:forwardeuler} with a smaller step size $\Delta \varepsilon$ (for a given value of $\Delta$ such that $\Delta \gg \varepsilon$). Then, we perform another iteration of \eqref{eq:forwardeuler} with step size $(1 - N\epsilon)\Delta$, where we assume that $N \ll 1/\epsilon$. The idea of this numerical scheme is to use the short step size iterations to allow the fast dynamics to converge to the equilibrium manifold $\mathcal{M}_0$, and then take a larger step size $(1 - N\epsilon)\Delta$ to let the slow subsystem evolve. 

During the simulation with the short step size $\Delta \varepsilon$, the discretized system is as follows:
\begin{equation}\label{eq:dtsmallstep}
    X_{k+1} = X_k + \Delta \epsilon F(X_k, u_k, \epsilon).
\end{equation}

After allowing the system to evolve according to \eqref{eq:dtsmallstep} for $N$ steps, the simulated system takes a larger step of length $(1 - N\varepsilon)\Delta$ according to
\begin{equation}\label{eq:dtbigstep}
    X_{k+N+1} = X_{k+N} + (1 - N\epsilon)\Delta F(X_{k + N}, u_{k + N}, \epsilon).
\end{equation}

Note that the total time interval covered by these $N+1$ steps is $\Delta$. This SMFE method is described in Algorithm \ref{alg:proposed_scheme}.
\begin{algorithm}[h]
\caption{Stabilized Multirate Forward Euler (SMFE) method for simulating the singularly perturbed system \eqref{eq:generalss}}
\label{alg:proposed_scheme}
\begin{algorithmic}[1]
\Require {$X(0)$, $F$, $u$, $\Delta$, $\epsilon$, $N$}
\For{$i = 0, 1, \dots$}
\For{$n = 1, \dots, N$}
\State $X(i\Delta + n\Delta \epsilon) \gets X(i \Delta + [n-1]\Delta \epsilon) +$ \\ \qquad\qquad\qquad\qquad\quad $\Delta \epsilon F(X(i \Delta + [n-1]\Delta \epsilon),u(i \Delta))$
\EndFor
\State $X([i+1] \Delta) \gets X(i \Delta + N \Delta \epsilon) +$ \\
\qquad\qquad\qquad\quad\quad $\Delta (1 - N \epsilon) F(X(i \Delta + N \Delta \epsilon),u(i \Delta))$
\EndFor
\end{algorithmic}
\end{algorithm}

\section{Stability analysis} \label{sec:stability_analysis}

In this section, we perform a stability analysis of the SMFE method and derive the requirements for the simulation step size $\Delta$ and the number of steps $N$. We demonstrate the efficiency of SMFE compared to the classic FEM while ensuring numerical stability.

We consider the block-triangular form in \eqref{eq:testtri}, where applying the SMFE method leads to the discrete-time system
\begin{align*}
\left[\begin{array}{l}
x(t + \Delta) \\
\eta(t + \Delta)
\end{array}\right] \approx \underbrace{\left( I + \Delta (1 - N \epsilon) A \right) \left( I + \Delta \epsilon A \right)^N}_{=: A_\Delta} \left[\begin{array}{l}
x(t) \\
\eta(t)
\end{array}\right].
\end{align*}

This discretized model is stable if and only if all the eigenvalues of $A_\Delta$ have magnitude less than one. Therefore, for stability, it is necessary and sufficient that the following condition holds for every eigenvalue $\lambda$ of $A$:
\begin{align} \label{eq:stability_condition}
|1 + \Delta(1 - N \epsilon) \lambda|\, |1 + \Delta \epsilon \lambda|^N < 1.
\end{align}

Given that $\epsilon \ll 1$ and that $A$ contains stable eigenvalues $\lambda_\text{slow} = O(1)$ and $\lambda_\text{fast} = O(1 / \epsilon)$, we first focus on the slowest eigenvalue $\hat{\lambda}_\text{slow}$, defined similarly to \eqref{eq:maxlamdafast}:
\begin{align}\label{eq:maxlamdaslow}
    \hat{\lambda}_\text{slow} = \mathop{\arg\min}\limits_{\lambda_{\text{slow}, i}}\{\mathrm{Re}(\lambda_{\text{slow}, i})\mid i = 1,2,\dots,n\}.
\end{align}
For small $\epsilon$ but fixed $\Delta$ and $N$, we use the approximation $(1 + \Delta \epsilon \lambda)^N \approx 1 + N \Delta \epsilon \lambda$ and neglected the $O(\epsilon)$ terms, thus condition~\eqref{eq:stability_condition} yields
\begin{align} \label{eq:approx_stab_condition_slow}
|1 + \Delta \hat{\lambda}_\text{slow}| < 1.
\end{align}
Therefore, for the slow subsystem, the stability condition for the Forward Euler discretization is recovered.

Next, we consider the fastest eigenvalue $\hat{\lambda}_\text{fast}$ as in \eqref{eq:maxlamdafast} and define the scaled eigenvalue $\tilde{\lambda} = \epsilon \hat{\lambda}_\text{fast} = O(1)$. For $\epsilon \ll 1$, condition~\eqref{eq:stability_condition} is equivalent to
\begin{align*}
N > - \frac{\ln |1 + \Delta \hat{\lambda}_\text{fast} - N \Delta \tilde{\lambda}|}{\ln |1 + \Delta \tilde{\lambda}|}.
\end{align*}

For very small $\epsilon$, while keeping $N$ and $\Delta$ fixed, it can be noticed that $\ln |1 + \Delta \hat{\lambda}_\text{fast} - N \Delta \tilde{\lambda}| > 0$, so for condition~\eqref{eq:stability_condition} to be satisfied we first need that $\ln |1 + \Delta \tilde{\lambda}| < 0$, \ie, $|1 + \Delta \tilde{\lambda}| < 1$. Under this assumption, considering $\epsilon \ll 1$, while keeping $\Delta$ fixed and assuming $N \ll 1/\epsilon$, we obtain the approximate condition
\begin{align} \label{eq:approx_stab_condition_fast}
N > - \frac{\ln (\Delta |\tilde{\lambda}| / \epsilon)}{\ln |1 + \Delta \tilde{\lambda}|} = O\left(\ln\left(\frac{1}{\epsilon}\right)\right).
\end{align}

Condition~\eqref{eq:approx_stab_condition_fast} confirms that the assumption $N \ll 1/\epsilon$ can be satisfied. Note that although the condition on $N$ is not independent of $\varepsilon$, its dependence is very mild due to the logarithm in \eqref{eq:approx_stab_condition_fast}.

\textbf{Remark 2.} \textit{Note that in condition~\eqref{eq:approx_stab_condition_fast}, as $\Delta \rightarrow \epsilon$, the lower bound on $N$ tends to dramatically increase, which shows that choosing a too small value for $\Delta$ can lead to a large number of simulation steps.}

\textbf{Remark 3.} \textit{In the case of a single fast eigenvalue, if one could choose $\Delta$ such that $1 + \Delta \tilde{\lambda} \approx 0$ while ensuring that $|1 + \Delta \lambda_\text{slow}| < 1$ for all slow poles, then the number of iterations $N$ can be reduced considerably. This is because the base in the second factor in \eqref{eq:stability_condition} becomes approximately zero.}

In terms of computational effort, since both $|1 + \Delta \hat{\lambda}_\text{slow}| < 1$ and $|1 + \Delta \tilde{\lambda}| < 1$ are required, the SMFE method requires
\begin{align*}
(N + 1) \frac{T}{\Delta}
= O\left( - \frac{\ln (\Delta |\tilde{\lambda}| / \epsilon)}{\ln |1 + \Delta \tilde{\lambda}|} \frac{T}{\Delta} \right)
=  O\left(T \hat{\lambda} \ln\left(\frac{1}{\epsilon}\right)\right)
\end{align*}
function evaluations in order to simulate the model for a simulation horizon $T$, since $\Delta = O(1/\hat{\lambda})$, where $\hat{\lambda}$ is the scaled eigenvalue of the largest magnitude of the system:
\begin{align*}
\hat{\lambda} = 
\begin{cases} 
\hat{\lambda}_\text{slow}, & \text{if } \operatorname{Re}(\hat{\lambda}_\text{slow}) < \operatorname{Re}(\epsilon \hat{\lambda}_\text{fast}), \\
\epsilon \hat{\lambda}_\text{fast}, & \text{otherwise}.
\end{cases}
\end{align*}

Notice again that this number depends only mildly on $\epsilon$. Therefore, SMES is shown to stabilize the explicit method for SPSs using far fewer iterations.

\section{Error analysis}\label{sec:error_analysis}

In this section, we analyze the numerical error characteristics of the SMFE method using the block-triangular form in \eqref{eq:testtri}. We aim to demonstrate that during the $N$ small-step simulation process, the system dynamics rapidly approach the equilibrium manifold of the corresponding reduced model. Consequently, the major simulation error originates from the discretization with the larger step $\Delta(1 - N\varepsilon)$.

First, we assume that the system in \eqref{eq:testtri} starts from $[x(t), \eta(t)] = [x_\Delta, \eta_\Delta]$.
After $N$ iterations with small step size $\Delta \epsilon$, the fast dynamics evolves as
\begin{align*}
  \eta(t + N\Delta \epsilon) = |1 + \Delta \epsilon \hat{\lambda}_\text{fast}|^N\eta_\Delta.
\end{align*}

Using the stability condition \eqref{eq:approx_stab_condition_fast} and substituting for $N$ yields
\begin{align}\label{eq:fastconverge}
    \eta(t + N\Delta \epsilon) & < |1 + \Delta \epsilon \hat{\lambda}_\text{fast}|^{- \frac{\ln (\Delta |\tilde{\lambda}| / \epsilon)}{\ln |1 + \Delta \tilde{\lambda}|}}\eta_\Delta.
\end{align}
Given $\tilde{\lambda} = \epsilon \hat{\lambda}_\text{fast}$ and the principle $a^b = e^{b \ln a}$ for positive $a$, \eqref{eq:fastconverge} can be simplified to
\begin{align}
    \eta(t + N\Delta \epsilon) &< \frac{\epsilon \eta_\Delta}{\Delta |\tilde{\lambda}|} = O(\epsilon).
\end{align}
While the fast dynamics converge rapidly to the order of $O(\epsilon)$, the slow dynamics evolve as
\begin{align*}
    x(t + N\Delta \epsilon) &= |1 + \Delta \epsilon \hat{\lambda}_\text{slow}|^N x_\Delta \\
    & \approx |1 + N \Delta \epsilon \hat{\lambda}_\text{slow}| x_\Delta \\
    & = x_\Delta + O(\epsilon),
\end{align*}
where we use the approximation $(1 + \delta)^N \approx 1 + N \delta$ for small $\delta$. This shows that the slow state remains close to $x_\Delta$ with changes of order $O(\varepsilon)$. Therefore, after $N$ small-step iterations, the state $[x(t + N\Delta \varepsilon), \eta(t + N\Delta \varepsilon)]$ converges to $[x_\Delta, 0]$ with an error of $O(\varepsilon)$. Note that $[x_\Delta, 0]$ lies on the equilibrium manifold $\mathcal{M}_0: \eta = 0$ of the discretized reduced model obtained by setting $\varepsilon = 0$ in \eqref{eq:testtri}. The discretized reduced model is given by
\begin{equation}\label{eq:testdisreduced}
    \begin{aligned}
        x(t+\Delta) &= |1 + \Delta(1 - N \epsilon) \lambda_\text{slow}|\, x(t + N\Delta \epsilon), \\
        \eta(t + \Delta) &= 0.
    \end{aligned}
\end{equation}
Note that \eqref{eq:testdisreduced} is a Forward Euler discretization with step length $\Delta(1 - N \epsilon)$ of the reduced model
\begin{equation}\label{eq:testreduced}
    \begin{aligned}
        &\dot{\Bar{x}} = ( A_{11} - A_{12}A^{-1}_{22}A_{21} )\, \Bar{x}(t), \\
        &\Bar{\eta}(t) = 0.
    \end{aligned}
\end{equation}
Since the $N$ small-step iterations introduce an error of $O(\varepsilon)$, the error from the Forward Euler discretization in \eqref{eq:testdisreduced}, which is $O(\Delta^2)$ \cite{hairer1991solving}, becomes the dominant source when $O(\Delta^2) \gg O(\varepsilon)$. Therefore, the error introduced by the SMFE method on the full model over time interval $\Delta$ is very close to the discretization error produced during the simulation of step length $\Delta(1 - N \varepsilon)$ on the reduced system. 

\section{Numerical illustration}\label{sec:numerical}

\begin{figure}
    \centering
    \includegraphics[width=0.75\linewidth]{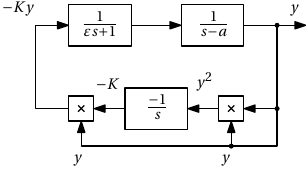}
    \caption{Adaptive control system with small parasitic dynamics \cite{kokotovic1999singular}.}
    \label{fig:adasys}
\end{figure}

In this section, we numerically illustrate the SMFE method on the simulation of a nonlinear adaptive control system with parasitic dynamics from \cite{kokotovic1999singular}, as shown in Fig.~\ref{fig:adasys}. The state-space model of the adaptive system is
\begin{equation}\label{eq:adaori}
\begin{aligned}
    \dot{y} & =a y+z, \\
\dot{k} & =y^2, \\
\varepsilon \dot{z} & =-z+u=-z-k y,
\end{aligned}
\end{equation}
where $a$ is the plant parameter and $\epsilon$ is a small positive parasitic time constant. A result of adaptive control theory \cite{narendra1980stable} is that the update law $\dot{k} =y^2$ for the adjustable gain $k$ stabilizes the system for all constant $a$ when the parasitic dynamics are neglected (\ie, $\epsilon = 0$). In this case, the adaptive system is reduced to the quasi-steady-state form
\begin{equation}\label{eq:adaslow}
\begin{aligned}
\dot{y} & =a y+z, \\
\dot{k} & =y^2, \\
z & = u = -k y,
\end{aligned}
\end{equation}
which can be numerically studied using explicit methods without resorting to very small time steps. However, for the original system in \eqref{eq:adaori}, explicit methods become inefficient due to numerical stability constraints, as discussed in Section~\ref{sec:FEM}. Nevertheless, since the parasitic dynamics are asymptotically stable for $0< \epsilon \ll 1$, they should converge rapidly to the equilibrium manifold $\mathcal{M}_0\colon z = -k y$. Therefore, we first apply the SMFE method to study the system with parasitic dynamics and compare the simulation results to those of the classical FEM to first address the efficiency of the SMFE method. The simulations are conducted with initial conditions $[y(0), k(0), z(0)] = [0, 0, 1]$, parameters $a = -1$ and $\epsilon = 10^{-6}$, over a simulation horizon $t \in [0, 5]$. 

\subsection{Efficiency of the SMFE method}

\begin{figure}
    \centering
    \includegraphics[width=1\linewidth]{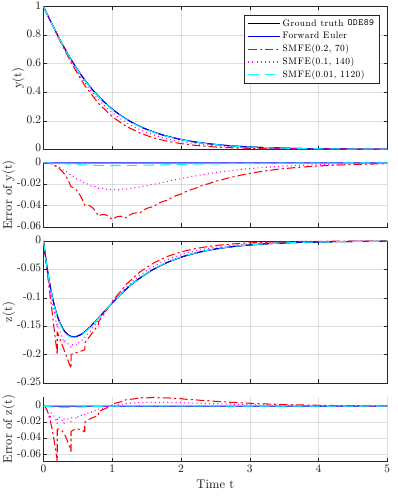}
    \caption{State trajectories of $y(t)$ and $z(t)$, and simulation errors with different methods. SMFE$(\Delta, N)$ represents the SMFE method using a large step of length $\Delta$ and small steps of number $N$.}
    \label{fig:adasimu}
\end{figure}

\begin{table*}[htbp]\label{table1}
\centering
\caption{Number of iterations and MSE for simulating the system in Fig.~\ref{fig:adasys} using different methods}
\begin{tblr}{
  row{odd} = {c},
  row{2} = {c},
  row{3} = {c},
  vline{2} = {-}{},
  hline{2} = {-}{},
  hline{3} = {-}{},
}
Method     & Forward Euler & SMFE(0.2, 70) & SMFE(0.2, 140) & SMFE(0.2, 1120) & SMFE(0.1, 140) & SMFE(0.1, 1120) & SMFE(0.01, 1120)\\
Iterations  & $5,000,000$ & $1,780$ & $3,530$ & $28,000$ & $7,050$ & $56,000$ & $561,000$ \\
MSE & $1.90\times10^{-14}$ & $8.29\times10^{-4}$ & $8.26\times10^{-4}$ & $8.25\times10^{-4}$ & $1.97\times10^{-4}$ & $1.96\times10^{-4}$ & $1.89\times10^{-6}$ \\
\end{tblr}
\end{table*}

The simulation results are shown in Fig.~\ref{fig:adasimu} and Table~\hyperref[table1]{\uppercase\expandafter{\romannumeral1}}. First, we simulate the system using \texttt{ODE89}, a high-order explicit Runge-Kutta method of order 9(8) that enhances the solution accuracy \cite{verner2010numerically}. We use the solution obtained by \texttt{ODE89} as the reference ``true'' solution of the system. In our simulations, we employ the SMFE method with different parameter settings, primarily adjusting the value of the large step size $\Delta$ and the number of small steps $N$, denoted as SMFE$(\Delta, N)$. Meanwhile, we use FEM with a step size small enough to guarantee its numerical stability. We consider the mean square error (MSE) of the simulated $y-z$ trajectories using different methods towards the ground truth as the error metric. 

From Fig.~\ref{fig:adasimu}, it can be observed that all four numerical methods provide stable simulations. However, as shown in Table~\hyperref[table1]{\uppercase\expandafter{\romannumeral1}}, achieving numerical stability with the classical FEM requires at least $5\times10^6$ iterations, whereas the SMFE method reduces the number of iterations to the order of $10^3$. This result confirms our stability analysis, demonstrating that the selection of the small step number $N$ depends only mildly on $\epsilon$ for numerical stability. Additionally, the significant reduction in the number of iterations does not lead to excessive error; the error can be further minimized by decreasing the large step size $\Delta$ and increasing the number of small steps $N$. 




\subsection{Error convergence}

\begin{figure*}[htbp]
    \centering
    \includegraphics[width=1\linewidth]{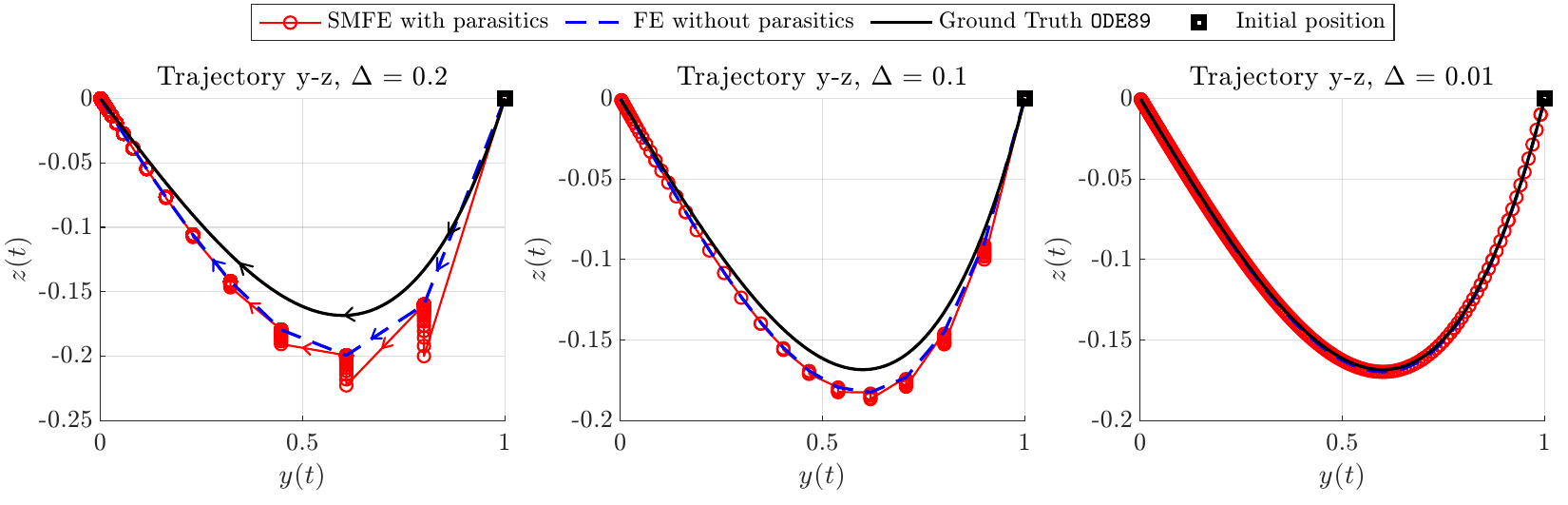}
    \caption{2D state trajectories of $y(t) - z(t)$ under different simulation schemes. SMFE is used to simulate the system with parasitic dynamics in \eqref{eq:adaori}, while FEM is used to simulate the reduced system in \eqref{eq:adaslow}. The large step size of SMFE and the step size of FEM are both $\Delta$ in each subplot. The direction of the state trajectories is indicated using arrows in the left subplot.}
    \label{fig:adaplot}
\end{figure*}

As shown in Table~\hyperref[table1]{\uppercase\expandafter{\romannumeral1}}, $\Delta$ and $N$ have different effects on the error of the SMFE method. When $\Delta$ is fixed, increasing $N$ only slightly reduces the error. However, when $N$ is fixed, reducing $\Delta$ significantly decreases the error. This observation is consistent with our analysis in Section~\ref{sec:error_analysis}: the main error of the SMFE method arises from the large-step discretization determined by $\Delta$, rather than the small-step discretization involving $N$ steps. Moreover, the error magnitude in Table~\hyperref[table1]{\uppercase\expandafter{\romannumeral1}} can be considered as of order $O(\Delta^2)$. 


To visualize the dependence of the error on $\Delta$, we plot the states $y(t)$ and $z(t)$ as 2D trajectories. We apply the SMFE method to the system with parasitic dynamics in \eqref{eq:adaori}, and then apply FEM on the reduced system in \eqref{eq:adaslow}. In each simulation, we set the step size of FEM equal to the large step size $\Delta$ used in the corresponding SMFE method. The system parameters and initial conditions remain the same as in the previous simulations. The simulation results are shown in Fig.~\ref{fig:adaplot}.

From Fig.~\ref{fig:adaplot}, we observe that as the step size $\Delta$ decreases, the trajectories obtained by both FEM and SMFE converge toward the true trajectory. Notably, in the subfigures (especially the first one), we see that while the large step sizes in the SMFE method cause the system trajectory to deviate from the equilibrium manifold $\mathcal{M}_0\colon z = -k y$ (represented by the blue trajectory), the subsequent small steps allow the system trajectory to quickly return to this manifold. This observation confirms the error analysis in Section~\ref{sec:error_analysis}: with the same step size, the numerical solution of SMFE approaches the equilibrium manifold $\mathcal{M}_0\colon z = -k y$ of the corresponding Forward Euler discretized reduced system.

In summary, the SMFE method provides stable numerical solutions for SPSs more efficiently than the classical FEM. While FEM can yield accurate simulation results, it necessarily requires a very small time step size due to its stability constraint, which may lead to excessively large simulation times. In contrast, the SMFE method allows users to balance the trade-off between error and simulation efficiency by adjusting $\Delta$ and $N$. Moreover, the MSE error remains within an acceptable range of order $O(\Delta^2)$.

\section{Conclusions}\label{sec:conclusion}

In this paper, we have addressed the challenge of numerically solving IVPs for SPSs by proposing SMES, which enhances classical explicit numerical methods through a multirate strategy with variable time steps, enabling efficient and stable simulations without the need for prohibitively small time steps or implicit formulations. Through stability and error analyses, we have demonstrated that SMES effectively stabilizes the fast dynamics, allowing the system to remain close to its equilibrium manifold. SMES provides a practical solution for efficiently simulating SPSs, offering a balance between computational cost and accuracy that classical explicit methods cannot achieve. Moreover, SMES is easy to implement and can be generalized to other explicit methods beyond the Forward Euler method used in this paper.


\bibliographystyle{ieeetr}
\bibliography{yiboshi}

\end{document}